\documentclass[12pt]{article}
\pdfoutput=1
\usepackage{amsmath}    
\usepackage{graphicx}   
\usepackage{verbatim}   
\usepackage{color}      
\usepackage{subfigure}  
\usepackage{amssymb}
\usepackage{cite}

\newtheorem{theorem}{Theorem}[section]
\newtheorem{proposition}[theorem]{Proposition}

\newtheorem{corollary}[theorem]{Corollary}
\newtheorem{lemma}[theorem]{Lemma}

\newtheorem{definition}[theorem]{Definition}

\newcommand{\proof}{\noindent{\bf Proof.\ }}
\newcommand{\qed}{\hfill $\square$\medskip}
\def\cp{\,\square\,}
\def\sp{\,\boxtimes\,}

\DeclareMathOperator {\gp} {gp}
\DeclareMathOperator {\ip} {ip}

\begin{document}
	
\title{The graph theory general position problem on some interconnection networks}
	
\author{
	Paul Manuel $^{a}$
	\and
	Sandi Klav\v zar $^{b,c,d}$
}

\date{}

\maketitle
\vspace{-0.8 cm}
\begin{center}
	$^a$ Department of Information Science, College of Computing Science and Engineering, Kuwait University, Kuwait \\
	{\tt pauldmanuel@gmail.com}\\
	\medskip
	
	$^b$ Faculty of Mathematics and Physics, University of Ljubljana, Slovenia\\
	{\tt sandi.klavzar@fmf.uni-lj.si}\\
	\medskip
	
	$^c$ Faculty of Natural Sciences and Mathematics, University of Maribor, Slovenia\\
	\medskip
	
	$^d$ Institute of Mathematics, Physics and Mechanics, Ljubljana, Slovenia\\
	
\end{center}

\begin{abstract}
Given a graph $G$, the (graph theory) general position problem is to find the maximum number of vertices such that no three vertices lie on a common geodesic. This graph invariant is called the general position number (gp-number for short)  of $G$ and denoted by ${\rm gp}(G)$. In this paper, the gp-number is determined for a large class of subgraphs of the infinite grid graph and for the infinite diagonal grid. To derive these results, we introduce monotone-geodesic labeling and prove a Monotone Geodesic Lemma that is in turn developed using the Erd\"os-Szekeres theorem on monotone sequences. The gp-number of the 3-dim infinite grid is bounded. Using isometric path covers, the gp-number is also determined for Bene\v{s} networks.
\end{abstract}


\medskip
\noindent{\bf AMS Subj.\ Class.: 05C12, 05C82} 

\section{Introduction}
\label{sec:introduction}

A set $S$ of vertices of a graph $G$ is called a {\em general position set} if no three vertices of $S$ lie on a common geodesic. A general position set $S$ of maximum cardinality is a {\em gp-set} of $G$ and its cardinality is the {\em general position number} (in short {\em gp-number}) of $G$ denoted by $\gp(G)$. The general position problem was introduced in~\cite{MaKl-2017+} and in particular motivated by the discrete geometry General Position Subset Selection Problem~\cite{Froese-2017+, PaWo13} which is to determine a largest subset of points in general position. The classical no-three-in-line problem however goes back all the way to Dudeney~\cite{dudeney-1917}; for more recent developments on it see~\cite{misiak-2016, PoWo07} and references therein.

In~\cite{MaKl-2017+}, several upper bounds on $\gp(G)$ were given. Connections between general position sets and packings were also investigated in order to obtain lower bounds on the gp-number. In addition, the general position problem was shown to be NP-complete. In this paper, we continue the study of the graph theory general position problem and focus on classes of interconnection networks. In order to determine their gp-number, a couple of new techniques are developed along the way. 

We proceed as follows. In the rest of this section definitions needed are listed. In the subsequent section, some results from~\cite{MaKl-2017+} are recalled. The concept of monotone-geodesic labellings is also introduced and a Monotone Geodesic Lemma is established. This lemma is derived from the Erd\"os-Szekeres theorem on monotone sequences. A couple of other techniques related to isometric subgraphs are also developed. Then, in Section~\ref{sec:networks}, the gp-number is determined for a large class of subgraphs of the grid graph (including the infinite grid itself) and for the infinite diagonal grid. A lower and an upper bound on the gp-number of the 3-dim grid is also given. In Section~\ref{sec:benes} the general position problem is solved for Bene\v{s} networks using isometric path covers. In the concluding section some directions for further study are suggested.  

Unless stated otherwise, graphs considered in this paper are connected. The {\em distance} $d_G(u,v)$ between vertices $u$ and $v$ of a graph $G$ is the number of edges on a shortest $u,v$-path. Shortest paths are also known as {\em geodesics} or {\em isometric paths}. A subgraph $H=(V(H),E(H))$ of a graph $G=(V(G),E(G))$ is {\em isometric} if $d_H(x,y) = d_G(x,y)$ holds for every pair of vertices $x,y$ of $H$. The size of a largest complete subgraph of a graph $G$ is its {\em clique number} $\omega(G)$. 

\section{Monotone-geodesic labeling}
\label{sec:monotone}

To approach the general position problem on interconnection networks, we first recall some known tools and then develop some new ones. First, the following simple fact will also be useful to us. 

\begin{proposition}
\label{prp:isometric}
Let $H$ be an isometric subgraph of a graph $G$. Then $S\subseteq V(H)$ is a general position set of $H$ if and only if $S$ is a general position set of $G$. 
\end{proposition}

\proof
Let $u,v,w\in V(H)$. Then $d_H(u,v) = d_H(u,w) + d_H(w,v)$ if and only if $d_G(u,v) = d_G(u,w) + d_G(w,v)$. That is, $u,v,w$ are on a common geodesic in $H$ if and only if they are on a common geodesic in $G$.  
\qed

An {\em isometric path cover} of a graph $G$ is a collection of geodesics that cover $V(G)$, cf.~\cite{Fitz99, pan-2006}. If $v$ is a vertex of a graph $G$, then let $\ip(v,G)$ be the minimum number of isometric paths, all of them starting at $v$, that cover $V(G)$. A vertex of a graph $G$ that lies in some gp-set of $G$ is called a {\em gp-vertex} of $G$. With these concepts in hand we can recall the following result. 

\begin{theorem} [\cite{MaKl-2017+}]
\label{thm:iso-v-path-gp}
If $R$ is a general position set of a graph $G$ and $v\in R$, then
\begin{equation}
\label{math:eq2}
|R| \leq \ip(v,G) + 1\,.
\end{equation}
In particular, if $v$ is a gp-vertex, then $\gp(G) \le \ip(v,G) + 1$. 
\end{theorem}

A sequence of real numbers is {\it monotone} if it is monotonically increasing or monotonically decreasing. The celebrated Erd\"os-Szekeres result, cf.~\cite[Theorem 1.1]{bukh-2014}, read as follows.  

\begin{theorem}[\cite{ErSz35}]
\label{thm:Er-Sze}
For every $n\ge 2$, every sequence $(a_1,\ldots, a_N)$ of real numbers with $N\ge (n-1)^2 + 1$ elements contains a monotone subsequence of length $n$. 
\end{theorem}

We will also say that a sequence $((x_1, y_1), \ldots (x_k, y_k))$ of points in the Cartesian plane is {\em monotone} if the sequences $(x_1, \ldots, x_k)$ and $(y_1, \ldots, y_k)$ are both monotone. For example $((1,4)$, $(2,4)$, $(5,3)$, $(5,2)$, $(6,1))$ is a monotone sequence. Theorem~\ref{thm:Er-Sze} has the following consequence tailored for us. 

\begin{corollary}
\label{cor:5-points}
If $n\in {\mathbb N}$ and $S$ is a set of $(n - 1)^2 + 1$ points in the Cartesian plane, then $S$ contains $n$ points that form a monotone sequence.   
\end{corollary}

\proof
Let $N = (n - 1)^2 + 1$ and let $S = \{(x_1,y_1),\ldots, (x_N,y_N)\}$ be an arbitrary set of $N$ points. We may assume without loss of generality that $x_1\le \cdots \le x_N$. By Theorem~\ref{thm:Er-Sze}, the sequence $(y_1, \ldots, y_N)$ contains a monotone subsequence of length $n$. This subsequence together with the corresponding first coordinates $x_i$ forms a required monotone sequence. 
\qed

If $n=3$, then Corollary~\ref{cor:5-points} asserts that any set of five points contains a monotone sequence of length $3$. For example, the set $\{(1,4), (2,3), (3,5), (3,2), (5,3)\}$ contains a monotone subsequence $((1,4), (2,3), (5,3))$. 

\begin{definition}[Monotone-geodesic labeling]
\label{Mon-Geo-Pro}
Let $G = (V(G),E(G))$ be a graph. Then an injective mapping $f:V(G)\rightarrow {\mathbb R}^2$ is a {\em monotone-geodesic labeling} of $G$ if the following holds: If $x$, $y$ and $z$ are vertices of $G$ such that the sequence of labels $(f(x)$, $f(y)$, $f(z))$ is monotone, then $x$, $y$, and $z$ lie on a common geodesic of $G$.  
\end{definition}

\begin{figure}[ht!]
	\begin{center}
		\scalebox{0.75}{\includegraphics{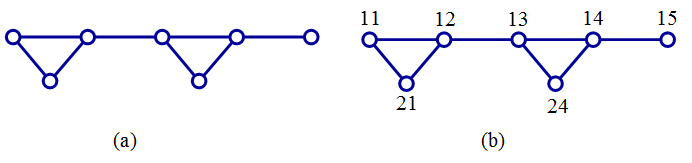}}
	\end{center}
	\caption{A graph equipped with a  monotone-geodesic labeling.}
	\label{fig:Monotone-Geodesic-Property}
\end{figure}

For an example see Fig.~\ref{fig:Monotone-Geodesic-Property}, where a graph is shown together with a monotone-geodesic labeling. 

We are now ready for the main insight of this preliminary section. 

\begin{lemma} [Monotone Geodesic Lemma]
\label{lem:TMGAM4}
If a graph $G$ admits a monotone-geodesic labeling, then $\gp(G) \le 4$. 
\end{lemma}
\proof
Suppose on the contrary that $S$ = $\{v_1,\ldots, v_5\}$ is a general position set of $G$. Let $f:V(G)\rightarrow {\mathbb R}^2$ be a monotone-geodesic labeling, where $f(v_i) = (x_i,y_i)$ for $i\in [5]$. Corollary~\ref{cor:5-points} applied for the case $n=3$ yields that $f(S)$ contains three points (labels) that form a monotone sequence, let it be $(f(v_{i_1}), f(v_{i_2}), f(v_{i_3}))$. Since $f$ is a monotone-geodesic labeling, the vertices  $v_{i_1}$, $v_{i_2}$, and $v_{i_3}$ lie on a common geodesic of $G$ which is a contradiction. 
\qed

From Lemma~\ref{lem:TMGAM4} it follows that not all graphs admit monotone-geodesic labelings. In particular, such a graph must necessarily have a small clique number.  

\begin{corollary}
\label{cor:clique}
If a graph $G$ admits a monotone-geodesic labeling, then $\omega(G)\le 4$. 
\end{corollary}

\proof
If $K$ is a complete subgraph of $G$, then $V(K)$ is (in view of Proposition~\ref{prp:isometric}) a general position set of $G$ and so $\gp(G)\ge \omega(G)$. Hence $\omega(G)\le \gp(G)\le 4$ by Lemma~\ref{lem:TMGAM4}.
\qed

Characterizing graphs that admit monotone-geodesic labellings seems to be an interesting open problem. It would also be interesting to characterize the graphs $G$ which satisfy $\omega(G) = \gp(G)$. 
 
\section{General position sets of grid networks}
\label{sec:networks}

By now we have prepared the main tools needed to determine (or bound) the gp-number of several interconnection networks that are based on the Cartesian and the strong product of graphs~\cite{hik-2011}. The {\em Cartesian product} $G\cp H$ of graphs $G$ and $H$ is the graph with the vertex set $V(G) \times V(H)$, vertices $(g,h)$ and $(g',h')$ being adjacent if either $g=g'$ and $hh'\in E(H)$, or $h=h'$ and $gg'\in E(G)$. The Cartesian product is a classical graph operation that is still intensively studied, cf.~\cite{bresar-2017, bog-2017, rall-2017, yang-2017}. The {\em strong product} $G\sp H$ is obtained from $G\cp H$ by adding, for every edge $gg'\in E(G)$ and every edge $hh'\in E(H)$, the edges $(g,h)(g',h')$ and $(g,h')(g',h)$. (We refer to~\cite{barrag-2016, zhao-2017} for a couple of recent developments on the strong product.) The {\em infinite $2$-dim grid} is the Cartesian product $P_{\infty}\cp P_{\infty}$ while the {\em infinite $2$-dim diagonal grid} is the strong product $P_{\infty}\sp P_{\infty}$. Using the standard notation from~\cite{hik-2011} we will denote them by $P_{\infty}^{\cp,2}$ and by $P_{\infty}^{\sp,2}$, respectively. Similarly, the {\em infinite $3$-dim grid} is the Cartesian product $P_{\infty}^{\cp,3}$.

\subsection{2-dim grids}

\begin{figure}[ht!]
	\begin{center}
	\scalebox{0.75}{\includegraphics{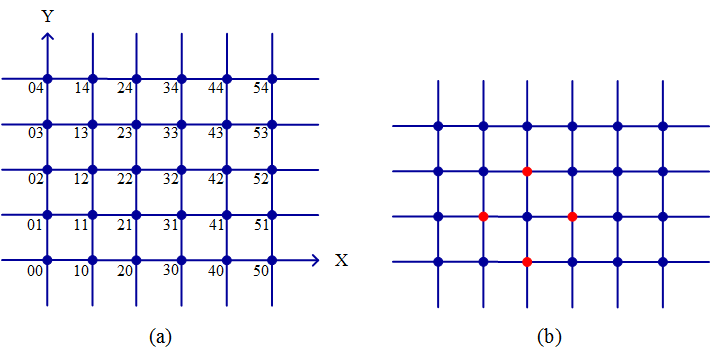}}
	\end{center}
	\caption{(a) The graph $P_{\infty}^{\cp,2}$ together with the natural labeling $f$ of its vertices, where $f(v_i,v_j) = (i,j)$ is briefly written as $ij$. (b) The red vertices form a general position set of $P_{\infty}^{\cp,2}$.}
	\label{fig:FIGR}
\end{figure}

Let $V(P_{\infty}) = \{\ldots, v_{-2}, v_{-1}, v_0, v_1, v_2, \ldots \}$ where $v_i$ is adjacent to $v_j$ if and only if $|i-j| = 1$. Then $V(P_{\infty}^{\cp,2}) = \{(v_i,v_j):\ i,j\in {\mathbb Z}\}$. Set now $f: V(P_{\infty}^{\cp,2}) \rightarrow {\mathbb R}^2$ with $f(v_i,v_j) = (i,j)$; see Fig.~\ref{fig:FIGR}(a). In this way the vertices of $P_{\infty}^{\cp,2}$ are labeled with the integer points in the Cartesian coordinate system. As this is a labeling of $P_{\infty}^{\cp,2}$ that (most probably) first comes to our minds, we call $f$ the {\em natural labeling} of $P_{\infty}^{\cp,2}$. 

A graph $G$ is a {\em grid graph} if it is an induced connected subgraph of $P_{\infty}^{\cp,2}$. Then we have: 

\begin{theorem}
\label{thm:special-grids}
Let $G$ be a grid graph equipped with the restriction $f|G$ of the natural labeling $f$ of $P_{\infty}^{\cp,2}$. If $G$ contains $P_3\cp P_3$ as a subgraph and $f|G$ is a monotone-geodesic labeling, then $\gp(G) = 4$. 
\end{theorem}

\proof
Since $f|G$ is a monotone-geodesic labeling, $\gp(G) \le 4$ by Lemma~\ref{lem:TMGAM4}. On the other hand a general position set of order $4$ as shown in Fig.~\ref{fig:FIGR}(b) exists in $G$ because it contains $P_3\cp P_3$ as a subgraph and since such a $P_3\cp P_3$ is necessarily an isometric subgraph of $G$.
\qed

Since the natural labeling $f$ of $P_{\infty}^{\cp,2}$ is monotone-geodesic, Theorem~\ref{thm:special-grids} yields: 

\begin{corollary}
\label{cor:TGP4IGR}
$\gp(P_{\infty}^{\cp,2}) = 4$. 
\end{corollary}

\subsection{2-dim diagonal grids}

\begin{figure}[ht!]
	\begin{center}
	\scalebox{0.76}{\includegraphics{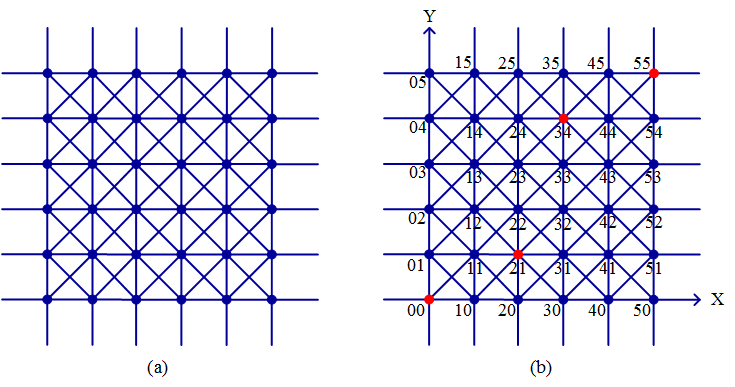}}
	\end{center}
	\caption{(a) The infinite 2-dim diagonal grid $P_{\infty}^{\sp,2}$. (b) The natural labeling of $P_{\infty}^{\sp,2}$ is not monotone-geodesic.}
	\label{fig:FIDGR1}
\end{figure}

We next consider the infinite 2-dim diagonal grid $P_{\infty}^{\sp,2}$, see Fig.~\ref{fig:FIDGR1}(a). One can label the vertices of $P_{\infty}^{\sp,2}$ with the natural labeling as used for  $P_{\infty}^{\cp,2}$, see Fig.~\ref{fig:FIDGR1}(b). However, now this natural labeling is no longer monotone-geodesic. For instance, the sequence $((0,0), (2,1), (3,4), (5,5))$ (see the red vertices in Fig.~\ref{fig:FIDGR1}(b)) is monotone (and so is every subsequence of it of length $3$), but no three of the corresponding vertices lie on a common geodesic. 
	
Despite the fact that the approach with the natural labeling does not work for $P_{\infty}^{\sp,2}$, we still have the following result. 

\begin{theorem}
\label{thm:TGP4IDGR}
$\gp(P_{\infty}^{\sp,2}) = 4$.
\end{theorem}

\begin{figure}[ht!]
	\begin{center}
	\scalebox{0.65}{\includegraphics{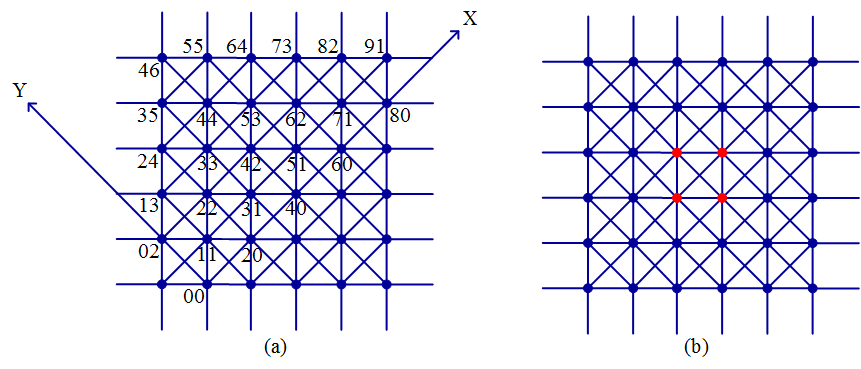}}
	\end{center}
	\caption{(a) A different labeling of $P_{\infty}^{\sp,2}$.(b) A general position set of $P_{\infty}^{\sp,2}$.}
	\label{fig:FIDGR2}
\end{figure}

\proof 
In order to show that $\gp(P_{\infty}^{\sp,2}) \le 4$, it is enough to identify a monotone-geodesic labeling for $P_{\infty}^{\sp,2}$. Consider the labeling of $P_{\infty}^{\sp,2}$ as shown in Fig.~\ref{fig:FIDGR2}(a) and call it $g$. Note that $g$ is derived from the natural labeling by rotating the Cartesian coordinate system by $45^{\circ}$.

We claim that the labeling $g$ is monotone-geodesic. So let $u,v,w$ be vertices of $P_{\infty}^{\sp,2}$ such that the sequence $(g(u), g(v), g(w))$ is monotone. We may assume without loss of generality that $g(u) = (0,0)$. Let $g(v) = (v_1, v_2)$ and  $g(w) = (w_1, w_2)$. Consider the case when $0\le v_1\le w_1$ and $0\le u_2\le w_2$. Then the vertex $v$ lies in the quadrant above the $x$ and $y$ coordinate axis (in the first quadrant), cf. Fig.~\ref{fig:FIDGR2}(a) again. Apply translation of axes from $u$ to $v$. Now $w$ lies in the first quadrant of the new translated coordinate system. It implies that $v$ lies on a $u,w$-geodesic. The other cases are symmetrical and can be argued similarly. Therefore, $\gp(P_{\infty}^{\sp,2}) \le 4$ by Lemma~\ref{lem:TMGAM4}. 

Since the red vertices from Fig.~\ref{fig:FIDGR2}(b) form a general position set, $\gp(P_{\infty}^{\sp,2}) \geq 4$. In conclusion, $\gp(P_{\infty}^{\sp,2}) = 4$.
\qed

\subsection{General position set of boron sheets}
\label{subsec:boron}

\begin{figure}[ht!]
	\begin{center}
		\scalebox{0.65}{\includegraphics{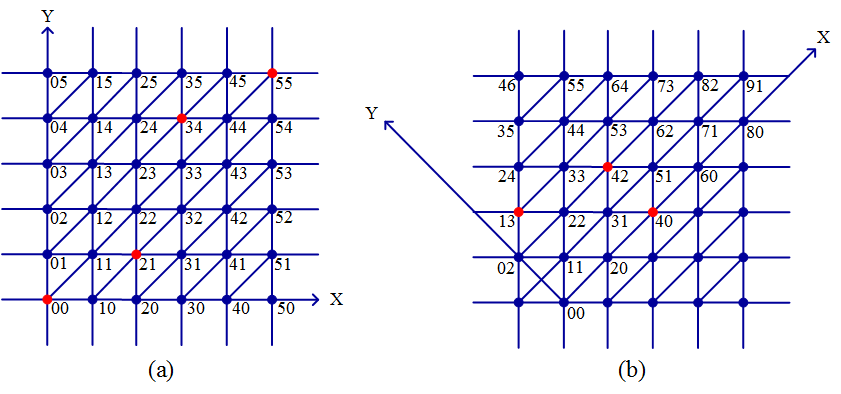}}
	\end{center}
	\caption{(a) Monotone-geodesic property fails  (b) Monotone-geodesic property fails again.}
	\label{fig:FBORSH}
\end{figure}

A boron sheet is a chemical graph which is also called a triangular grid. An infinite boron sheet is described in Fig.~\ref{fig:FBORSH}. The 2-dim grid  $P_{\infty}^{\cp,2}$ is a subgraph of the infinite triangular grid which is in turn a subgraph of the 2-dim diagonal grid $P_{\infty}^{\sp,2}$. Since the gp-number of the $P_{\infty}^{\cp,2}$ and $P_{\infty}^{\sp,2}$ is 4, it may lead to conclude that the gp-number of the triangular grid is also 4.

In the previous subsections we have seen two different types of labeling of grids. The conventional type of labeling for boron sheets is illustrated in Fig.~\ref{fig:FBORSH}(a) and the second type is illustrated in Fig.~\ref{fig:FBORSH}(b). Surprisingly, monotone-geodesic property fails in either case. For example,  the sequence $((0,0), (2,1), (3,4), (5,5))$ in Fig.~\ref{fig:FBORSH}(a) is monotone, yet no three of the corresponding vertices lie on a common geodesic. In the same way, the sequence $((1,3), (4,2), (4,0))$ from Fig.~\ref{fig:FBORSH}(b) is monotone, but the there vertices do not lie on a geodesic.

\begin{figure}[ht!]
	\begin{center}
		\scalebox{1.0}{\includegraphics{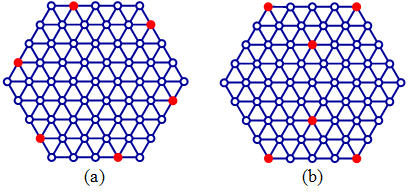}}
	\end{center}
	\caption{(a) A general position set of order $6$ (b) Another example of a general position set of the same cardinality.}
	\label{fig:FHexBor6}
\end{figure}

Thus, the gp-number of infinite boron sheet may be more than $4$. In fact, we can show that the gp-number of infinite boron sheet is at least $6$, see Fig.~\ref{fig:FHexBor6}.

The carbon nano sheets which are also called hexagonal or honeycomb grids behave in the same way as boron sheets. We believe that the gp-number of boron sheets and carbon nano sheets is $6$. 

\subsection{3-dim grids}
\label{subsec:3grids}

We next consider the infinite $3$-dim grid, that is, the graph $P_{\infty}^{\cp,3}$. In view of Corollary~\ref{cor:TGP4IGR} one might expect that either $\gp(P_{\infty}^{\cp,3}) = 2\cdot 3 = 6$ or $\gp(P_{\infty}^{\cp,3}) = 2^3 = 8$.\footnote{These were actually the guesses of the audience in the University Newcastle, Australia, when one of the authors was presenting the results of this paper.} Hence the main result of the subsection comes as a surprise. 

We start with the following simple yet useful result to be applied in identifying general position sets in $P_{\infty}^{\cp,3}$.  

\begin{lemma}
	\label{lem:LGPCOND}
	Let $G=(V(G),E(G))$ be a graph and $S\subseteq V(G)$. If there exists an integer $k$ such that $k \leq d(x,y) < 2k$ holds for every different $x,y\in S$, then $S$ is a general position set.  
\end{lemma}

\proof
Suppose $S$ is not general position set. Then there exist vertices $x$, $y$, and $z$ of $S$ such that $d(x,y) = d(x,z) + d(z,y)$. Since $d(x,z)\ge k$ and $d(z,y)\ge k$ we have $d(x,y)\geq 2k$, a contradiction to the lemma's hypothesis. 
\qed

Lemma~\ref{lem:LGPCOND} for $k=1$ says that the vertex set of any complete subgraph of a graph forms a general position set. Note also that if diameter of $G$ is at most $3$, then Lemma~\ref{lem:LGPCOND} (for $k=2$) asserts that every independent set of $G$ is a general position set.  

Now we are ready for the announced surprising result. 

\begin{proposition}
\label{prp:T3DGRID}
$10\le \gp(P_{\infty}^{\cp,3}) \le 16$. 
\end{proposition} 

\proof 
For the lower bound it suffices to construct a general position set of order $10$. Consider $P_{5}^{\cp,3}$ equipped with the natural labeling of its vertices set $S = \{(2,2,0)$, $(3,1,1)$, $(1,3,1)$, $(2,0,2)$, $(0,2,2)$, $(4,2,2)$, $(2,4,2)$, $(1,1,3)$, $(3,3,3)$, $(2,2,4)\}$. Note that here (and in the rest of the proof) we have identified the vertices with the points in $3$-dim Euclidean space. Now, it is easy to verify that $3 \leq d(x,y) \leq 5$ for every pair of vertices  $x, y \in S$. Thus, by Lemma~\ref{lem:LGPCOND}, $S$ is a general position set. Since $P_{5}^{\cp,3}$ is an isometric subgraph of $P_{\infty}^{\cp,3}$, Proposition~\ref{prp:isometric} implies that $S$ is a general position set of $P_{\infty}^{\cp,3}$. 

For the upper bound consider an arbitrary set $S = \{(x_i,y_i,z_i):\ i\in [17]\}$ of vertices of $P_{\infty}^{\cp,3}$ or order $17$. We may without loss of generality assume that $x_1\le x_2\le\cdots \le x_{17}$. By Theorem~\ref{thm:Er-Sze}, the sequence $(y_1, y_2, \ldots, y_{17})$ contains a monotone subsequence of order $5$, say $(y_{i_1}, \ldots, y_{i_5})$. Using Theorem~\ref{thm:Er-Sze} again, the sequence $(z_{i_1}, \ldots, z_{i_5})$ contains a monotone subsequence, say $(z_{i_{j_1}}, z_{i_{j_2}}, z_{i_{j_3}})$. But now, the vertices $(x_{i_{j_1}}, y_{i_{j_1}}, z_{i_{j_1}})$, $(x_{i_{j_2}}, y_{i_{j_2}}, z_{i_{j_2}})$, and $(x_{i_{j_3}},y_{i_{j_3}},z_{i_{j_3}})$ lie on a geodesic, so $S$ is not a general position set.
\qed

Inductively using the argument from the second part of the proof of Proposition~\ref{prp:T3DGRID} we can also infer the following result:

\begin{proposition}
	\label{prp:TKDGRID} If $k$ is an arbitrary positive integer, then 
	$\gp(P_{\infty}^{\cp,k}) < \infty\,.$ 
\end{proposition} 

In addition to finding the exact value of $\gp(P_{\infty}^{\cp,3})$, the gp-problem for $P_{\infty}^{\sp,3}$ is also worth-studying. Needless to mention that the gp-problem of $\gp(P_{\infty}^{\cp,k})$ and  $\gp(P_{\infty}^{\sp,k})$, where $k\ge 3$, will remain a challenge to researchers.

\section{General position sets in products of cycles}
\label{sec:torus}

\begin{figure}[ht!]
	\begin{center}
		\scalebox{0.8}{\includegraphics{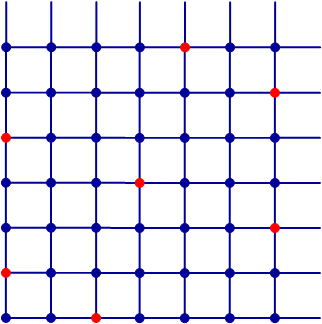}}
	\end{center}
	\caption{A general position set with seven vertices for torus.}
	\label{fig:FTORGPS}
\end{figure}

A {\em torus} is a Cartesian product of two cycles, cf.\ Fig.~\ref{fig:FTORGPS} where $C_7\cp C_7$ is illustrated such that the ``long edges" are only indicated. 

\begin{lemma}
	\label{lem:LGPS3}
	Let $S$ be general position set of $P_r\cp P_s$. If $S$ contains a corner vertex (that is, a vertex of degree $2$), then $|S| \leq 3$. 
\end{lemma}

\begin{figure}[ht!]
	\begin{center}
		\scalebox{0.8}{\includegraphics{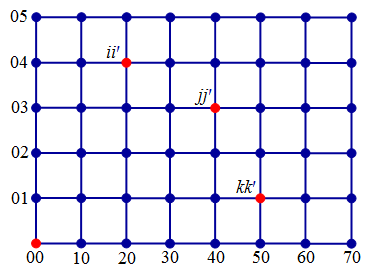}}
	\end{center}
	\caption{Suppose $(i,i')$, $(j,j')$, and $(k,k')$ are the other three vertices of $S$ where $i \leq j \leq k$.}
	\label{fig:FLGPS3}
\end{figure}

\proof
Suppose on the contrary that $|S|$ = 4. Assume without loss of generality that $(0,0) \in S$. Let $(i,i')$, $(j,j')$, and $(k,k')$ be the other three vertices of $S$ where $i \leq j \leq k$.  
If $i = j$, then $(0,0)$, $(i,i')$, and $(j,j')$ lie in a geodesic. So $i < j$. Similarly, we see that $j < k$. Now, if $i' \leq j'$, then $(0,0)$, $(i,i')$, and $(j,j')$ would again be in a geodesic, therefore $i' > j'$. Similarly, we also see that $j' > k'$; see Fig.~\ref{fig:FLGPS3}. 

But then $(i,i')$, $(j,j')$, and $(k,k')$ lie in a geodesic, a final contradiction.
\qed

\begin{theorem}
	\label{TTOR79}
$7\le \gp(C_r\cp C_s)\le 9$. 
\end{theorem}

\begin{figure}[ht!]
	\begin{center}
		\scalebox{0.8}{\includegraphics{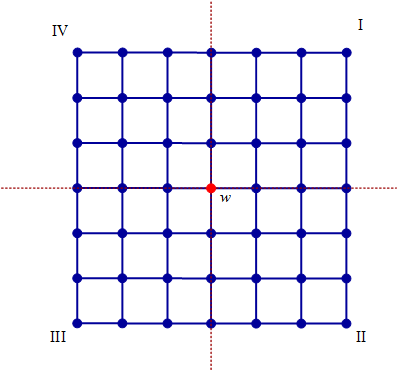}}
	\end{center}
	\caption{Four quadrants of a torus with respect to $w$.}
	\label{fig:FTORPART}
\end{figure}

\proof 
Let $S$ be general position set of $C_r\cp C_s$. Let $w$ be an arbitrary vertex from $S$. Partition $C_r\cp C_s$ into $4$ parts with respect to $w$ as illustrated in Fig.~\ref{fig:FTORPART}. 

By Lemma~\ref{lem:LGPS3}, in each of the four parts the set $S$ can have a maximum of two vertices other than $w$. Thus the gp-number of $C_r\cp C_s$ is at most $9$. ON the other hand, Fig.~\ref{fig:FTORGPS} demonstrates the existence of a general position set of order $7$. 

\qed

\section{Bene\v{s} networks}
\label{sec:benes}

In this section we determine the gp-number of Bene\v{s} networks. These networks are significant among inter-connection networks because they are rearrangeable non-blocking networks. (A network is rearrangeable non-blocking if any permutation can be realized by edge-disjoint paths when the entire permutation is known.) 

\begin{figure}[ht!]
	\begin{center}
		\scalebox{0.45}{\includegraphics{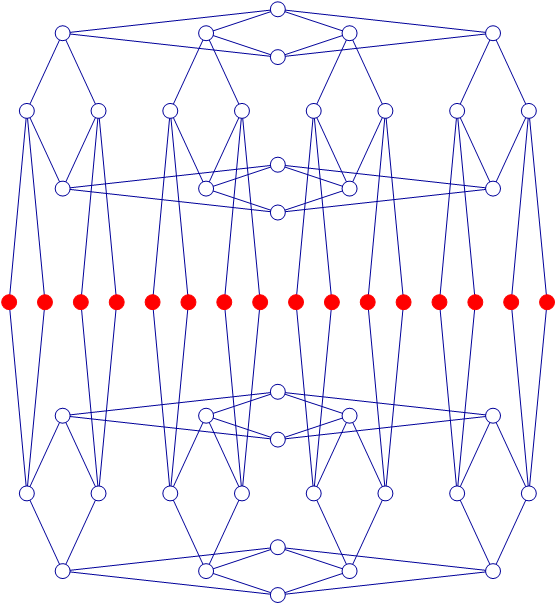}}
	\end{center}
	\caption{The $3$-dim Bene\v{s} network $BN(3)$. Its $2$-degree vertices are marked in red color and form a gp-set.}
	\label{fig:Benes-3-dim}
\end{figure}

The Bene\v{s} networks consist of back-to-back butterflies~\cite{MaAb08}, where in turn the $r$-dim butterfly has $n=2^{r}(r+1)$ nodes arranged in $r+1$ levels of $2^{r}$ nodes each. Each node has a distinct label $\langle w,i\rangle $, where $i$ is the level of the node $(1\leq i\leq r+1)$ and $w$ is a $r$-bit binary number that denotes the column of the node. Two nodes $\langle w,i\rangle $ and $\langle w',i'\rangle $ are linked by an edge if $i'=i+1$ and either $w$ and $w'$ are identical or $w$ and $w'$ differ only in the bit in position $i'$. We refer to~\cite[Section 11.4]{xu-2013} for basic properties of butterfly networks and to~\cite{burckel-2014, KlMa-2016} for a recent application and the average distance of these networks, respectively. Now, for $r\ge 1$ the $r$-dim Bene\v{s} network $BN(r)$ is constructed by merging two $r$-dim butterfly networks as shown in Fig.~\ref{fig:Benes-3-dim} for the case $r=3$. 

\begin{theorem}
\label{thm:Benes}
If $r\ge 1$, then $\gp(BN(r)) = 2^{r+1}$.
\end{theorem}

\proof
The case $r=1$ can be easily verified directly. In the rest let $r\ge 2$, let $R$ be an arbitrary general position set of $BN(r)$, and let $S$ be the set of degree $2$ vertices of $BN(r)$. See Fig.~\ref{fig:Benes-3-dim}, where the vertices of $S$ are drawn in red color. 

We will inductively show that $\gp(BN(r)) \le 2^{r+1}$ and for this sake, we distinguish two cases. 

\begin{figure}[ht!]
	\begin{center}
		\scalebox{0.45}{\includegraphics{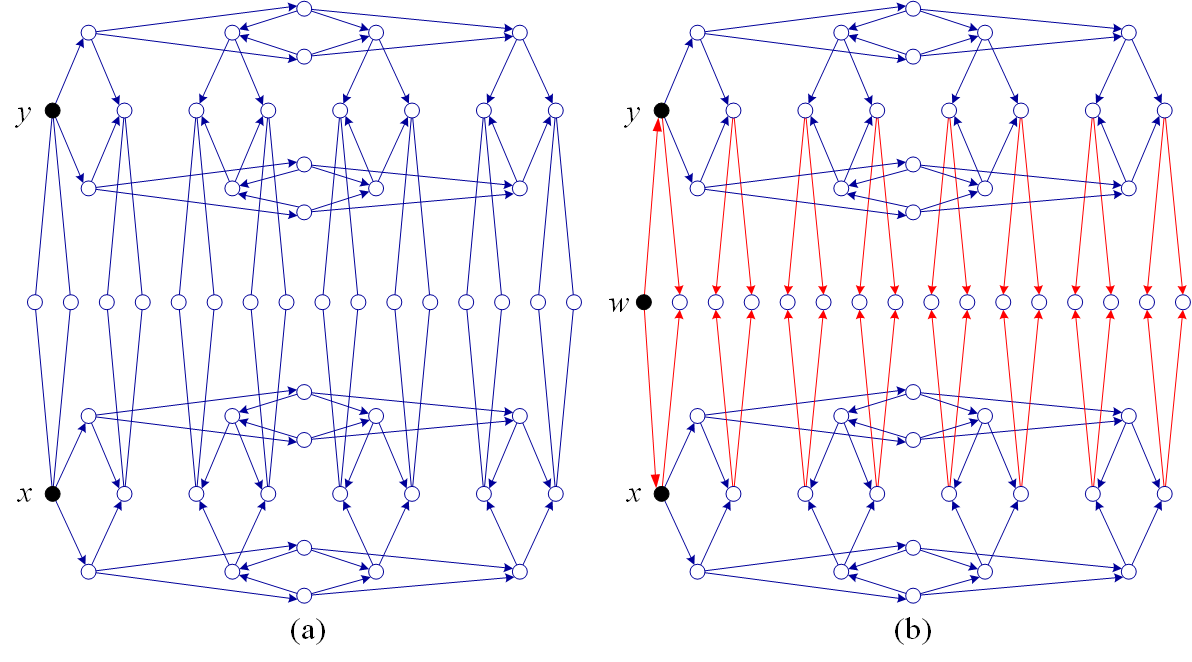}}
	\end{center}
	\caption{(a) Construction of $\Psi_x$ and $\Psi_y$ at inductive step $k = r-1$. (b) Construction of $\Psi_w$ at inductive step $k = r$.}
	\label{fig:Benes-Psi_w}
\end{figure}

\medskip\noindent
{\bf Case 1}: $R\cap S \neq \emptyset$. \\
Let $w\in R\cap S$. Then we inductively construct an isometric path cover 
$$\Psi_w = \{P_{wv}:\ v \in S, v \neq w, P_{wv}\ \mbox{is\ a\ fixed}\ w,v\mbox{-geodesic}\}$$ 
as follows. Let $x$ and $y$ be the two vertices of $BN(r)$ adjacent to $w$. Removing $S$ from $BN(r)$ leaves two $(r-1)$-dim Bene\v{s} networks $BN(r-1)$. By induction hypothesis, we can construct isometric path covers $\Psi_x$ and $\Psi_y$ of $BN(r-1)$, see Fig.~\ref{fig:Benes-Psi_w}(a). Then extend $\Psi_x$ and $\Psi_y$ to construct $\Psi_w$ of $BN(r)$, see Fig.~\ref{fig:Benes-Psi_w}(b).

Since $\Psi_w$ is an isometric path cover of $BN(r)$ and $w\in R$, Theorem~\ref{thm:iso-v-path-gp} implies that 
$$|R| \le \ip(w,BN(r)) + 1 \le |\Psi_w| + 1 = |S| = 2^{r+1}\,.$$

\medskip\noindent
{\bf Case 2}: $R \cap S = \emptyset$. \\
In this case, no vertex of $R$ has degree $2$ in $BN(r)$. Removing all the vertices of $S$ from $BN(r)$, the graph $BN(r)$ is disconnected into two $(r-1)$-dim Bene\v{s} networks $BN(r-1)$. By induction hypothesis, $\gp(BN(r-1)) \le 2^{r}$. Since the two copies of $BN(r-1)$ are isometric subgraphs of $BN(r)$, Proposition~\ref{prp:isometric} implies that the restriction of $R$ to each of the copies of $BN(r-1)$ contains at most $2^{r}$ vertices. Therefore, $|R| \le 2^{r+1}$. 

We have thus proved that $\gp(BN(r)) \le 2^{r+1}$. On the other hand, the set $S$ is a general position set of $BN(r)$ and we are done because $|\mathcal{S}| = 2^{r+1}$.
\qed

\section{Further research}
\label{sec:Fur-Research}
One of the key concepts of this paper is the monotone-geodesic labeling. A characterization of graphs that admit monotone-geodesic labelings will be very useful not only for the general position problem but also for other related topics. We have established a tool to test whether a given vertex set is a general position set. Using this result, it is demonstrated that the gp-number of the infinite $3$-dim grid is between $10$ and $16$. However, the exact $gp$-number of 3-dim grids is still unknown. As it is pointed out in  Subsection~\ref{subsec:3grids}, the gp-problem of $\gp(P_{\infty}^{\cp,k})$ and  $\gp(P_{\infty}^{\sp,k})$ will remain a challenge to researchers.
	
The general position problem for Bene\v{s} networks is solved using isometric path covers.  A Bene\v{s} network is a back-to-back butterfly network. However, the strategy applicable to  Bene\v{s} networks does not work for butterfly networks. It remains as a challenge to prove that the $gp$-number of $r$-dim butterfly is $2^r$.

The general position problem for 2-dim grids and 2-dim diagonal grids is solved using monotone-geodesic labellings and Monotone Geodesic Lemma. The structure of triangular grids (also called boron sheets, see~\cite{feng-2016, penev-2012}) is between $2$-dim grids and $2$-dim diagonal grids. The natural intuition is that the gp-number of triangular grids is $4$, because the gp-number of 2-dim grids and 2-dim diagonal grids is 4. However, the gp-number of the triangular grids is at least 6 and we conjecture that it is actually equal to $6$. 

\section*{Acknowledgment}
	
This work was supported and funded by Kuwait University, Kuwait. 


\end{document}